\documentclass[letterpaper]{article}
\usepackage{fullpage}
\usepackage{pslatex}
\usepackage{color}
\usepackage{graphics}
\usepackage{eepic}
\usepackage{amssymb}

\newcommand{\cmt}[1]{\ifhmode\newline\fi{\sf *** \ \ #1 \\}}

\newtheorem{theorem}{Theorem}

\newtheorem{definition}{Definition}
\newtheorem{lemma}{Lemma}
\newtheorem{corollary}{Corollary}
\newtheorem{remark}{Remark}
\newtheorem{fact}{Fact}

\newtheorem{example}{Example}
\newcommand{\qed}{\hspace{7pt}\rule{4pt}{5pt}}
\newlength{\saveparindent}
\def\bproof{\begin{rm}\protect\vspace{5pt}\noindent\textbf{Proof: }%
\addtolength{\parskip}{4pt}\setlength{\parindent}{0pt}}
\def\eproof{\end{rm}\addtolength{\parskip}{-4pt}%
\setlength{\parindent}{\saveparindent}}

\newcommand{\bprooff}[1]{\begin{rm}\protect\vspace{5pt}%
\noindent\textbf{Proof of #1: }\addtolength{\parskip}{4pt}%
\setlength{\parindent}{0pt}}

\newenvironment{proof}{\par\bproof}{\eproof\(\qed\) \par\medskip}

\newlength{\saveparskip}
\newcommand{\setof}[1]{\{#1\}}
\newcommand{\ZZ}{\mathbb{Z}}

\newcommand{\sgn}[1]{\mathrm{sign}(#1)}

\newcommand{\set}[2]{\left\{#1 \mbox{ : } #2 \right\}}
\newcommand{\smidge}{{\kern .05em}}

\newcommand{\barG}{\overline{G}}

\newcommand{\barU}{\overline{U}}
\newcommand{\barV}{\overline{V}}
\newcommand{\baru}{\overline{u}}
\newcommand{\barv}{\overline{v}}
\newcommand{\Toep}{\mathbb{T}}
\newcommand{\shape}[1]{\mathbf{shp}\left(#1\right)}

\begin{document}

\title{Distribution of the Size of a Largest Planar Matching and
Largest Planar Subgraph in Random Bipartite Graphs
%\thanks{}
}

\author{
{\sc Marcos Kiwi}\thanks{Gratefully acknowledges the support of
        MIDEPLAN via ICM-P01--05, and CONICYT via FONDECYT 1010689 
        and FONDAP in Applied Mathematics.}
        \\
          {\footnotesize Depto.~Ing.~Matem\'{a}tica and}\\[-1.5mm]
          {\footnotesize Ctr.~Modelamiento Matem\'atico UMI 2807, }\\[-1.5mm]
          {\footnotesize University of Chile}\\[-1.5mm]
          {\footnotesize Correo 3, Santiago 170--3, Chile}
        \\[-1.5mm]{\footnotesize e-mail: \texttt{mkiwi@dim.uchile.cl}}
\and
{\sc Martin Loebl}\thanks{Gratefully acknowledges the support of ICM-P01-05.
This work was done while visiting the Depto.~Ing.~Matem\'atica, U. Chile.}
        \\
           {\footnotesize Dept.~of Applied Mathematics and}\\[-1.5mm]
           {\footnotesize Institute of Theoretical Computer Science (ITI)}\\[-1.5mm]
           {\footnotesize  Charles University}\\[-1.5mm]
           {\footnotesize  Malostransk\'{e} n\'{a}m. 25, 118~00~~Praha~1}\\[-1.5mm]
           {\footnotesize Czech Republic}
        \\[-1.5mm]{\footnotesize e-mail: \texttt{loebl@kam.mff.cuni.cz}}
}

\date{}

\maketitle 

\begin{abstract}
We address the following question: When a 
  randomly chosen regular bipartite multi--graph is drawn in the plane
  in the ``standard way'', what is the distribution of its maximum
  size planar matching (set of non--crossing disjoint edges)
and maximum size planar subgraph (set of non--crossing edges
which may share endpoints)?
The problem is a generalization of the Longest Increasing Sequence (LIS)
  problem (also called Ulam's problem).
We present combinatorial identities which relate the number of $r$-regular
  bipartite multi--graphs with maximum planar matching (maximum planar subgraph)
of at most $d$ edges
  to a signed sum of restricted lattice walks in $\ZZ^d$, and to the number of 
  pairs of standard Young tableaux of the same shape and with a 
  ``descend--type'' property. 
Our results are obtained via generalizations 
  of two combinatorial proofs through which Gessel's identity can 
  be obtained (an identity that is crucial in the derivation of a 
  bivariate generating function associated to the distribution of 
  LISs, and key to the analytic attack on Ulam's problem).
\end{abstract} 


\section{Introduction}\label{sec.intro}
Let $U$ and $V$ henceforth denote two disjoint totally
 ordered sets (both ordered relations will be referred to by $\preceq$).
Typically, we will consider the case where $|U|=|V|=n$ and  
  denote the elements of $U$ and $V$ by $u_1,u_2,\ldots,u_n$
  and $v_1,v_2,\ldots,v_n$ respectively.
Henceforth, we will always assume that the latter enumeration 
  respects the ordered relation in $U$ or $V$, i.e.,
  $u_1\preceq u_2\preceq \ldots \preceq u_n$ and 
  $v_1\preceq v_2\preceq \ldots \preceq v_n$.

Let $G=(U,V;E)$ denote a bipartite multi--graph with color classes $U$ and $V$.
Two distinct edges $uv$ and $u'v'$ of $G$ are said to be \emph{noncrossing}
  if $u$ and $u'$ are in the same order as $v$ and $v'$; in other
  words, if $u\prec u'$ and $v\prec v'$ or $u'\prec u$ and $v'\prec v$.
A matching of $G$ is called \emph{planar} if every distinct pair
  of its edges is noncrossing.
We let $L(G)$ denote the number of edges of a maximum
 size (largest) planar matching in $G$ (note that $L(G)$ depends on
 the graph $G$ \emph{and} on the ordering of its color classes).
 
For the sake of simplicity we will concentrate solely in the case
  where $|E|=rn$ and $G$ is $r$--regular.

\medskip
When $r=1$, an $r$--regular multi--graph with color classes $U$ and $V$
  uniquely determines a permutation.
A planar matching corresponds thus to an increasing sequence of the
  permutation, where an increasing sequence
  of length $L$ of a permutation $\pi$ of $\{1,\ldots,n\}$ is
  a sequence 
  $1\leq i_1 < i_2 < \ldots < i_L \leq n$ such that 
  $\pi(i_1) < \pi(i_2) < \ldots < \pi(i_L)$.
The Longest Increasing Sequence (LIS)
  problem concerns the determination of the asymptotic, on $n$,
  behavior of the LIS for a randomly and uniformly chosen 
  permutation $\pi$.
The LIS problem is also referred to as ``Ulam's problem'' 
  (e.g., in~\cite{kingman73,bdj99,okounkov00}).
Ulam is often credited for raising it in~\cite{ulam61}
  where he mentions (without reference) 
  a ``well--known theorem'' asserting that given
  $n^{2}+1$ integers in any order, it is always possible to find among 
  them a monotone subsequence of $n+1$
  (the theorem is due to  Erd\H{o}s and Szekeres~\cite{es35}).
Monte Carlo simulations are reported in~\cite{bb67}, where it is 
  observed that over the range $n\leq 100$, the limit of the LIS
  of $n^2+1$ randomly chosen elements, 
  when normalized by $n$, approaches $2$.
Hammersley~\cite{hammersley72} gave a rigorous proof of the existence 
  of the limit and conjectured it was equal to $2$.
Later, Logan and Shepp~\cite{ls77}, based on a result by
  Schensted~\cite{schensted61}, proved that $\gamma\geq 2$;
  finally, Vershik and Kerov~\cite{vk77} obtained that $\gamma\leq 2$.
In a major recent breakthrough due to Baik, Deift, Johansson~\cite{bdj99}
  the asymptotic distribution of the LIS has been determined.
For a detailed account of these results, history and related work
  see the surveys of Aldous and Diaconis~\cite{ad99} and
  Stanley~\cite{stanley02}.

{From} the previous discussion, it follows that one way 
  of generalizing Ulam's problem is to study  
  the distribution of the size of the largest planar matching in 
  randomly chosen $r$--regular bipartite multi--graphs
  (for a different generalization see~\cite{steele77,bw88}).
This line of research, 
  originating in~\cite{kl02}, turns out to be relevant for
  the study of several other issues
  like the Longest Common Subsequence problem (see~\cite{klm05}),
  interacting particle systems~\cite{seppalainen98}, 
  digital boiling~\cite{gtw01}, 
  and is directly related to topics such as 
  percolation theory~\cite{alexander94} and random matrix 
  theory~\cite{johansson99}.

\subsection{Main Results}
We establish combinatorial identities which express $g(n;d)$ --- 
  the number of $r$-regular bipartite multi--graphs with planar 
  matchings with at most $d$ edges --- 
  in terms of:
\begin{itemize}
\item The number of 
  pairs of standard Young tableaux of the same shape and 
  with a ``descend-type'' property (Theorem~\ref{thm:perm-young}).

\item A signed sum of restricted lattice walks in $\ZZ^d$ 
  (Theorem~\ref{thm.1}).

% MARTIN: The paragraph can not go as an item of the above list.
%         Read the phrase before the list!! It says that what
%         follows are expressions for g(n;d)
%
%         I re-wrote the following paragrah AFTER the the list.
%         See below.
%\item
%We show that our arguments can be extended in order to
%  characterize the distribution of the largest size of 
%  planar subgraphs of randomly chosen $r$--regular bipartite multi--graphs
%(Theorem~\ref{thm.plk}).
\end{itemize}
Our arguments can be extended in order to
  characterize the distribution of the largest size of 
  planar subgraphs of randomly chosen $r$--regular bipartite multi--graphs
  (Theorem~\ref{thm.plk}).

\subsection{Models of Random Graphs: From k-regular Multi--graphs to Permutations}
Most work on random regular graphs is based on the so called random
  configuration model of Bender and Canfield and
  Bollob\'as~\cite[Ch.~II, \S~4]{bollobas85}.  
Below we follow this approach, but first we need to adapt the 
  configuration model to the bipartite graph scenario.  
Given $U$, $V$, $n$ and $r$ as above, let
  $\barU=U\times [r]$ and $\barV=V\times [r]$.
An $r$--configuration of $U$ and $V$ is 
  a one--to--one pairing of $\barU$ and $\barV$.  
These $rn$ pairs are called edges of the configuration.  
Hence, a configuration can be considered a graph, specifically,
  a perfect matching with color classes $\barU$ and $\barV$.
Moreover, viewing a configuration as such bipartite graph enables
  us to speak also about its planar matchings (here the total
  ordering on $\barU=U\times [r]$ and $\barV=V\times [r]$ is the lexicographic
  one induced by $\preceq$ and $\leq$).

The natural projection of $\barU=U \times [r]$ and $\barV=V\times [r]$ onto  
  $U$ and $V$ respectively (ignoring the second coordinate) projects 
  each configuration $F$ to a bipartite multi--graph $\pi(F)$ with 
  color classes $U$ and $V$.
Note in particular that $\pi(F)$ may contain multiple edges 
  (arising from sets of two or more edges in $F$ whose end--points 
  correspond to the same pair of vertex in $U$ and $V$).  
However, the projection of the uniform distribution over
  configurations of $U$ and $V$ is not the uniform
  distribution over all $r$--regular bipartite multi--graphs on $U$ and
  $V$ (the probability of obtaining a given multi--graph is proportional
  to a weight consisting of the product of a factor $1/j!$ for each
  multiple edge of multiplicity $j$).
Since a configuration $F$ can be considered a graph, it makes perfect
  sense to speak of the size $L(F)$ of its largest planar matching.

We denote an element $(u,i)\in \barU$ by $u^i$ and
  adopt an analogous convention for the elements of $\barV$.
We shall further abuse notation and denote by $\preceq$ the total
  order on $\barU$ given by $u^i \preceq \tilde{u}^j$ if
  $u\prec \tilde{u}$ or $u=\tilde{u}$ and $i\leq j$.
We adopt a similar convention for $\barV$.

Let $G_{r}(U,V;d)$ denote the set of all $r$--regular bipartite
  multi--graphs on $U$ and $V$ whose largest planar matching is of size
  at most $d$. 
Note that if $|U|=|V|=n$, then the cardinality of $G_{r}(U,V;d)$ depends 
  on $U$ and $V$ solely through $n$.
Thus, for $|U|=|V|=n$, let $g(n;d)=|G_{r}(U,V;d)|$.

The first step in our considerations is an identification of $G_{r}(U,V;d)$
with a subset of configurations of $U$ and~$V$.
Specifically, we associate to an $r$--regular multi--graph $G=(U,V;E)$ 
the $r$--configuration $\barG$ of $U$ and $V$ such that $\pi(\barG)=G$ where: 
If $(u,v)$ is an edge of multiplicity $t$ in $G$ for which there are $i$ edges
  $(u,v')$ in $G$ such that $v\prec v'$, 
  and $j$ edges $(u',v)$ in $G$ such that 
  $u\prec u'$, then for every $s\in [t]$, the pairing 
  $(u^{i+s},v^{j+t-s+1})$ belongs to $\barG$.
Note that the number of edges of $\barG$ equals the number of edges 
  of $G$.

Let $\barG_r(U,V;d)$ be the collection of configurations $\barG$ 
  associated to some $G\in G_{r}(U,V;d)$. Observe, that $g(n;d)=|\barG_{r}(U,V;d)|$.

For an edge $(\baru,\barv)$ 
  we say that $M\subseteq \set{\baru'\in \barU}{\baru'\preceq \baru}
    \times \set{\barv'\in \barV}{\barv'\preceq \barv}$
  is a planar matching that ends with $(\baru,\barv)$ if the edges in $M$ 
  are non--crossing and $(\baru,\barv)\in M$. 
Since there is a unique edge incident to every node in $\barG$, say
  $(\baru,\barv)$, we speak of a largest planar matching of $\barG$
  up to $\baru$ (or $\barv$) in order to refer to a largest planar matching
  that ends with edge $(\baru,\barv)$.

Note that the way in which $\barG$ is derived from $G$, implies
  in particular that for $u\in U$ and 
  $i\leq j$, the size of the maximum planar matching in $\barG$ 
  using nodes up to $u^i$ is at least as 
  large as the size of the maximum planar matching using
  nodes up to $u^j$.
A similar fact holds for elements $v\in V$.

\begin{figure}
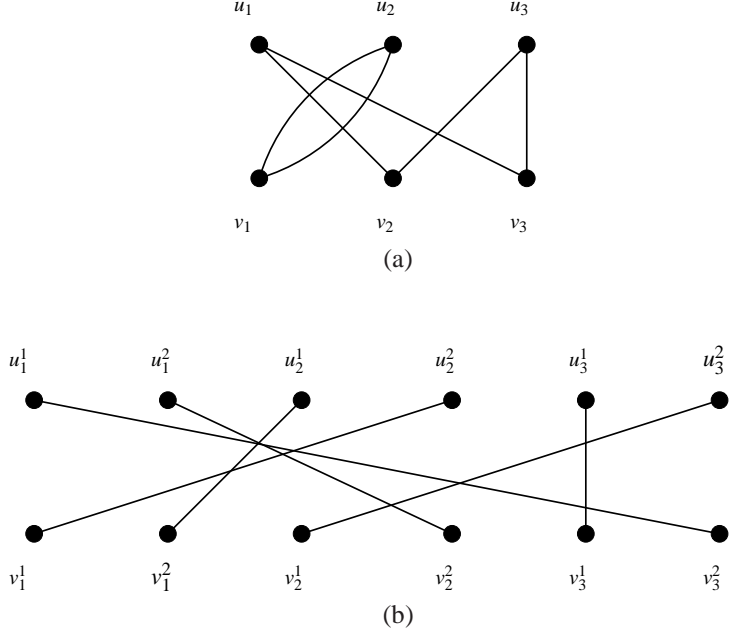

\begin{center}
\input graph.pstex_t \\
(a) 
\end{center}

\vspace{8pt}
\begin{center}
\input configuration.pstex_t \\
(b)
\end{center}
\caption{(a) A $2$--regular multi--graph $G$. (b) Configuration 
  $\barG$ associated to $G$.}
\label{fig:example-graphs}
\end{figure}

\medskip
Several of the concepts introduced in this section are illustrated
  in Figure~\ref{fig:example-graphs}.

\subsection{Young tableaux}\label{sec:young}
A \emph{(standard) Young tableau of shape} 
  $\lambda=(\lambda_1,\ldots,\lambda_r)$
  where $\lambda_1\geq \lambda_2\geq \ldots \geq \lambda_r\geq 0$,
  is an arrangement $T=(T_{k,l})$ 
  of $\lambda_1+\ldots+\lambda_r$ distinct integers
  in an array of left--justified rows, with $\lambda_i$ elements in 
  row $i$, such that the entries in each row are in increasing order from 
  left to right, and the entries of each column are increasing from 
  top to bottom (here we follow the usual convention that considers
  row $i$ to be above row $i+1$).
One says that $T$ has $r$ rows and $c$ columns 
  if $\lambda_r>0$ and $c=\lambda_1$ respectively.
The shape of $T$ will be henceforth denoted $\shape{T}$ and 
  the collection of Young tableau with entries in the set $S$
  and with at most $d$ columns will be denoted $T(S;d)$.

\begin{figure}
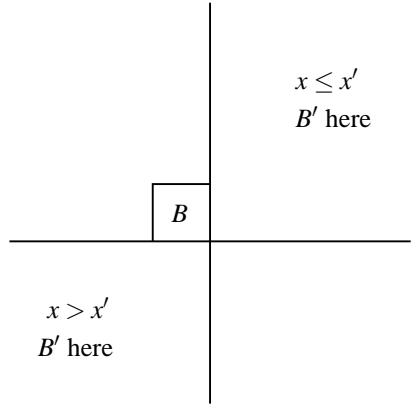

\begin{center}
\input young.pstex_t
\end{center}
\caption{New tableau entries created through row--insertions.}
\label{fig:rsk}
\end{figure}

The Robinson correspondence (rediscovered independently by Schensted)
  states that the set of permutations of $[m]$ is in one to one
  correspondence with the collection of pairs of equal shape tableaux 
  with entries in $[m]$.
The correspondence can be constructed through the Robinson--Schensted--Knuth
  (RSK) algorithm --- also referred to as row--insertion or row--bumping 
  algorithm.
The algorithm takes a tableau $T$ and a positive integer $x$, and constructs
  a new tableau, denoted $T\leftarrow x$. 
This tableau will have one more box than $T$, and its entries will be 
  those of $T$ together with one more entry labeled $x$, but there is 
  some moving around, the details of which are not of direct concern
  to us, except for the following fact:
\begin{lemma}\textbf{[Bumping Lemma~\cite[pag.~9]{fulton97}]}\label{lem:rsk}
Consider two successive row--insertions, first row inserting $x$ in a 
  tableau $T$ and then row--inserting $x'$ in the resulting tableau
  $T\leftarrow x$, given rise to two new boxes $B$ and $B'$ as shown in
  Figure~\ref{fig:rsk}.
\begin{itemize}
\item If $x\leq x'$, then $B$ is strictly left of and weakly below $B'$.
\item If $x> x'$, then $B'$ is weakly left of and strictly below $B$.
\end{itemize}
\end{lemma}
Given a permutation $\pi$ of $[m]$, the Robinson--Schensted--Knuth (RSK)
  correspondence constructs $(P(\pi),Q(\pi))$ such that
  $\shape{P(\pi)}=\shape{Q(\pi)}$ by,
\begin{itemize}
\item starting with a pair of empty
  tableaux, repeatedly row--inserting the elements
  $\pi(1),\ldots,\pi(n)$ to create $P(\pi)$, and,
\item placing the value $i$ into the box of $Q(\pi)$'s diagram 
  corresponding to the box created during the $i$--th insertion 
  into $P(\pi)$.
\end{itemize}
Two remarkable facts about the RSK algorithm which we will exploit
  are:
\begin{remark}\textbf{[RSK Correspondence~\cite[pag.~40]{fulton97}]}%
\label{rem:correspondence}
The RSK correspondence sets up a one--to--one mapping between 
  permutations of $[m]$ and pairs of tableaux $(P,Q)$ with
  the same shape.
\end{remark}

\begin{remark}\textbf{[Symmetry Theorem~\cite[pag.~40]{fulton97}]}%
\label{rem:transpose}
If $\pi$ is a permutation of $[m]$, then $P(\pi^{-1})=Q(\pi)$ and 
  $Q(\pi^{-1})=P(\pi)$.
\end{remark}
Moreover, it is easy to see that the following holds:
\begin{remark}\label{rem:ascending}
Let $\pi$ be a permutation of $[m]$.
Then, $\pi$ has no ascending sequence of 
  length greater than $d$ if and only if $P(\pi)$ and $Q(\pi)$ have 
  at most $d$ columns.
\end{remark}

The reader interested on an in depth discussion of Young tableaux is
  referred to~\cite{fulton97}.

% MARTIN: I moved the rest of this section to the beginning of 
%         the section ``First Proof''. It makes more sense, since
%         this way there are no proofs in the introduction. Moreover,
%         now each subsection 1.2, 1.3 and 1.4 concerns only one
%         topic and not their inter-relation.
%

\subsection{Walks}
We say that $w=w_0\ldots w_m$ is a lattice walk in $\ZZ^d$ of length $m$ if
  $||w_{i}-w_{i-1}||_1=1$ for all $1\leq i\leq m$.
Moreover, we say that $w$ starts at the origin and ends in $\vec{p}$ if
  $w_0=\vec{0}$ and $w_m=\vec{p}$.
For the rest of this paper, all walks are to be understood as lattice
  walks in $\ZZ^d$.
Let $W(d,m;\vec{p})$ denote the set of all walks of length $m$ from
  the origin to $\vec{p}\in\ZZ^d$.

We will often identify the walk $w=w_0\cdots w_m$ 
  with the sequence 
  $d_1\ldots d_m$ such that $w_{i}-w_{i-1}=\sgn{d_i}\vec{e}_{|d_i|}$, 
  where $\vec{e}_j$ denotes the $j$--th element of the canonical basis
  of~$\ZZ^d$.
If $d_i$ is negative, then we say that the $i$--th step is a 
  negative step in direction $|d_i|$, or negative step for short.
We adopt a similar convention when $d_i$ is positive.

We say that two walks are equivalent if both subsequences of the
  positive and the negative steps are the same.  
For each equivalence class consider the representative for 
  which the positive steps precede the negative steps. 
Each such representative walk may hence be written as 
  $a_{1}a_{2}\cdots |b_{1}b_{2}\cdots$ where
  the $a_i$'s and $b_j$'s are all positive.
For an arbitrary collection of walks $W$, all 
  with the same number of positive and the same number of negative
  steps, we henceforth denote by $W^*$ the collection of the 
  representative walks in $W$.

\medskip
Recall that one can associate to a permutation $\pi$ of $[d]$  
  the Toeplitz point $T(\pi)= (1-\pi(1),\dots,d-\pi(d))$.  
Note that in a walk from the origin to a Toeplitz point, the number of
  steps in a positive direction equals the number of steps in a
  negative direction.  
In particular, each such walk has an even length.

\medskip
In cases where we introduce notation for referring to 
  a family of walks from the origin 
  to a given lattice point $\vec{p}$, such as $W(d,m;\vec{p})$, 
  we sometimes consider instead of $\vec{p}$ a subset of lattice points $P$.
It is to be understood that we are thus making reference to
  the set of all walks in the family that end at 
  a point in $P$.
A set of lattice points of particular interest to the ensuing
  discussion is the set of Toeplitz points, henceforth denoted
  $\Toep$.

We now come to a simple but crucial observation: there is a 
  natural identification of $U\times [r]$ with $[rn]$ that 
  respects the total order in each of these sets ($\preceq$ in the
  former and $\leq$ in the latter).
A similar observation holds for $V\times [r]$.
Hence, when $m=rn$ the sequences of positive and negative steps in 
  a walk in $W(d,2m;\Toep)$ can be referred to as:
\[
a_{u_1^1}\cdots a_{u_1^r}a_{u_2^1}\cdots a_{u_2^r}
       \cdots a_{u_n^1}\cdots a_{u_n^r}
\quad\mbox{ and }\quad
b_{v_1^1}\cdots b_{v_1^r}b_{v_2^1}\cdots b_{v_2^r}
       \cdots b_{v_n^1}\cdots b_{v_n^r}\,.
\]
Let $W'(d,2m;T(\pi))$ be the set of all walks in $W^*(d,2m;T(\pi))$ 
  whose positive steps $a_{u_1^1}\cdots a_{u_n^r}$ and negative steps 
  $b_{v_1^1}\cdots b_{v_n^r}$ satisfy:
  $a_{u^i} \geq a_{u^{i+1}}$ and $b_{v^i} \geq b_{v^{i+1}}$ for all $u\in U$,
  $v\in V$ and $1\leq i< r$.

\section{Counting Planar Mathchings and Planar Subgraphs}
We are now ready to state the main result of this paper.
\begin{theorem}\label{thm.1}
\[
g(n;d)= \sum_{\pi}\sgn{\pi} \left|W'(d,2rn;T(\pi))\right|\,.
\]
\end{theorem}
%\begin{corollary}\label{cor.1}
%\[
%{{2rn}\choose{rn}} g(n;d) = 
%  \sum_{\pi}\sgn{\pi} \left|W^{*}(d,2rn;T(\pi))\right|\,.
%\]
%\end{corollary}
Our proof of Theorem~\ref{thm.1} is strongly based on the arguments
  used in~\cite{gww98} to prove the following result concerning
  $1$--regular bipartite graphs:
\begin{theorem}\label{thm.mot}
The signed sum of the number of walks 
  of length $2m$ from the origin to Toeplitz
  points is ${2m}\choose{m}$ times the number $u_m(d)$ of permutations 
  of length $m$ that have no increasing sequence
  of length bigger than $d$.
\end{theorem}
This last theorem gives a combinatorial proof of the following
  well known result:
\begin{theorem}{[Gessel's Identity]}\label{thm.mt2}
If $I_{\nu}(t)$ denotes the Bessel function of imaginary argument, then
\[
\sum_{m\geq 0}\frac{u_m(d)}{(m!)^2}x^{2m}=
\det(I_{|r-s|}(2x))_{r,s= 1,\dots, d}.
\]
\end{theorem}
% MARTIN: I like a lot your idea of moving here the section I added.
%         But not the way it was done. So I merged what was at 
%         the end of the section ``Walks'' sith the section ``Planar
%         Subgraph'' and created a this new section ``Counting ...''
%
%\section{Planar Subgraphs}
%
We now describe a random process which 
  researchers have studied, either explicitly or implicitly,
  in several different contexts. 
Let $X_{i,j}$ be a non--negative random
  variable associated to the lattice point $(i,j)\in [n]^{2}$. 
For $C\subseteq [n]^2$, we referr to $\sum_{(i,j)\in C} X_{i,j}$ 
  as the weight of $C$.
We are interested on the determination of the 
  distribution of the maximum weight of $C$
  over all $C=\setof{(i_1,j_1),(i_2,j_2),\ldots }$
  such that $i_1,i_2,\ldots$ and $j_1,j_2,\ldots$ are strictly increasing.

Johansson~\cite{johansson99} considered the case where the $X_{i,j}$s are 
  independent identically distributed according to a geometric 
  distribution.
Sep\"al\"ainen~\cite{seppalainen98} and Gravner, Tracy and Widom~\cite{gtw01}  
  studied the case where the $X_{i,j}$s are 
  independent identically distributed Bernoulli random variables
  (but, in the latter paper, the collections of lattice points 
  $C=\setof{(i_1,j_1),(i_2,j_2),\ldots }$ were such that 
  $i_1,i_2,\ldots$ and 
  $j_1,j_2,\ldots$ were weakly and strictly increasing respectively). 

The main result of this paper, i.e., Theorem~\ref{thm.1}, says that if
  $(X_{i,j})_{(i,j)\in [n]^2}$ is uniformly distributed over all adjacency
  matrices of $r$--regular multi--graphs, 
  then the distribution of the maximum weight evaluated at $d$ 
  can be expressed as a signed sum of 
  restricted lattice walks in $\ZZ^d$.
A natural question is whether a similar result holds if one relaxes
  the requirement that the sequences 
  $i_1,i_2,\ldots $ and $j_1,j_2,\ldots$ are strictly increasing.
For example, if one allows them to be weakly increasing.
% MARTIN: I modified the following phrase in order introduce the
%         notion of planar subgraph. It was not defined before,
%         but was being used in the paper's title and in a 
%         section title.
This is equivalent to 
  asking for the distribution of the size of a planar subgraph,
  i.e., the largest set of 
  non--crossing edges which may share endpoints in
  a uniformly chosen $r$--regular multigraph.
A line of argument similar to the one we will use in the 
  derivation of Theorem~\ref{thm.1} yields:
\begin{theorem}
\label{thm.plk}
Let 
  $\hat{g}(n;d)$ be the number of $r$-regular bipartite multi--graphs 
  with no larger than $d$ set of non--crossing edges which may share 
  endpoints.
Then, $\hat{g}(n;d)$ equals the number of pairs of equal shape
  Young tableaux in $T([rn];d)$ satysifying:
\begin{quote}
\textbf{Condition (\^{T}):} 
  If for each $i\in [n]$ and $1\leq s <r$, the row 
  containing $r(i-1)+s+1$ is weakly above the row containing
  $r(i-1)+s$.
\end{quote}
Moreover,
\[
\hat{g}(n;d) = \sum_{\pi}\sgn{\pi} \left|\widehat{W}'(d,2rn;T(\pi))\right|\,,
\]
where $\widehat{W}'(d,2rn;T(\pi))$ is the set of all walks in 
  $W^{*}(d,2rn;T(\pi))$ whose positive steps 
  $a_{u_{1}^1} \cdots a_{u_{n}^r}$ and negatives steps
  $b_{v_{1}^1} \cdots b_{v_{n}^r}$ satisfy:
  $a_{u^i} < a_{u^{i+1}}$ and $b_{v^i} < b_{v^{i+1}}$ for all $u\in U$,
  $v\in V$ and $1\leq i< r$.
\end{theorem}
 
In the rest of the paper we give two independent proofs 
  of Theorem~\ref{thm.1}.

\section{First Proof}\label{sec:first}
%
% MARTIN: The start of this section was what used to be at the 
%         end of the section on Young tableaux
%
Let $m=rn$.
Recall that $\barG_{r}(U,V;d)$ can be thought of as a collection
  of permutations of $[m]$. 
Thus, we may think of the RSK correspondence as being defined 
  over $\barG_{r}(U,V;d)$.
In particular, for an $r$--configuration $F$ of $U$ and $V$ we may
  write $(P(F),Q(F))$ to denote the pair of Young tableaux associated
  to the permutation determined by $F$.
Figure~\ref{fig:perm-young} shows the result of applying the RSK 
  algorithm to an $r$--configuration.
\begin{figure}
\begin{center}
\input perm-young.pstex_t
\end{center}
\caption{Pair of Young tableaux associated through the RSK algorithm 
  to the $2$--configuration of Figure~\ref{fig:example-graphs}.b.}
\label{fig:perm-young}
\end{figure}

We say that a Young tableau in $T([m];d)$ satisfies 
\begin{quote}
\textbf{Condition (T):} If 
  for each $i\in [n]$ and $1\leq s <r$, the row 
  containing $r(i-1)+s$ is strictly above the row containing
  $r(i-1)+s+1$.
\end{quote}
The following result characterizes the image of $\barG_{r}(U,V;d)$
  through the RSK correspondence.

\begin{theorem}\label{thm:perm-young}
The number $g(n;d)$ equals the number of pairs of equal shape tableaux 
  in $T([m];d)$ satisfying condition (T). 
Specifically, the RSK correspondence establishes a one--to--one 
  correspondance between $\barG_{r}(U,V;d)$ and the collection 
  of pairs of equal shape tableaux in $T([m];d)$ satisfying 
  condition~(T).
\end{theorem}

\begin{proof}
Let $\barG\in \barG_r(U,V;d))$.
Remark~\ref{rem:ascending} implies that 
  $P(\barG)$ and $Q(\barG)$ are tableaux of equal shape that
  belong to $T([m];d)$.
Corollary~\ref{lem:rsk} implies that, for every $i\in[n]$
  the row insertion process through which $P(\barG)$ is built is 
  such that the insertion of the values $r(i-1)+1,\ldots, ri$ 
  gives rise to a sequence of boxes each of which is strictly below
  the previous one.
This implies that $Q(\barG)$ satisfies condition (T).

We still need to show that $P(\barG)$ also satisfies condition (T).
For $\barG\in G_{r}(U,V;d)$ let the transpose of $\barG$, 
  denoted $\barG^{T}$, be the bipartite graph over color classes $U$ and 
  $V$ such that $u_iv_j$ is an edge of $\barG^{T}$ if and only if
  $u_jv_i$ is an edge of $\barG$.
Note that $\barG^T\in\barG_{r}(U,V;d)$ if and only if
  $\barG\in\barG_{r}(U,V;d)$.
A direct consequence of Remark~\ref{rem:transpose} is that 
  $(P(\barG^{T}),Q(\barG^{T}))=(Q(\barG),P(\barG))$.
Hence, $P(\barG)$ must also satisfy condition (T).

Suppose now that $(P,Q)$ is a pair of equal shape 
  tableaux in $T([m];d)$ both of which 
  satisfy condition~(T).
Let $F$ be an $r$--configuration of $U$ and $V$ such that 
  $(P(F),Q(F))=(P,Q)$
  (here we identify $u_{i}^s$ and $v_i^s$ and view $F$ as a permutation
  of $[m]$).
The existence of $F$ is guaranteed by Remark~\ref{rem:correspondence}. 
Remark~\ref{rem:ascending} implies that $F$'s largest planar
  matching is of size at most $d$.
Moreover, since $Q(F)$ satisfies property (T), Lemma~\ref{lem:rsk} 
  implies that the edges of $F$ incident to $u_{i}^{s}$
  and $u_{i}^{s+1}$ cross.
Similarly, one can conclude that the edges fo $F$ incident to
  $v_{i}^{s}$ and $v_{i}^{s+1}$ cross.
It follows that $F$ belongs to $\barG_r(U,V;d)$.
\end{proof}
\begin{example}
Note that condition (T), as guaranteed by Theorem~\ref{thm:perm-young}, is 
  reflected in the tableaux shown in Figure~\ref{fig:perm-young} 
  (for the tableau in the left; $4$, $1$ and $3$ are strictly above 
  $6$, $2$ and $5$ respectively, while for the tableau in the 
  right; $1$, $3$ and $5$ are strictly above $2$, $4$ and $6$ respectively).
\end{example}

For a walk $w=a_{1}\cdots a_{m}|b_{1}\cdots b_{m}$ in 
  $W'(d,2m;T(\pi))$ let 
  $\tilde{w}= \tilde{a}_{1}\cdots \tilde{a}_{m}|
        \tilde{b}_{1}\cdots \tilde{b}_{m}$
  be such that $\tilde{a}_{i}=a_{i}$ and 
  $\tilde{b}_{i}=b_{m-i}$.
Denote by $\widetilde{W}(d,2m;T(\pi))$ the collection of all $\tilde{w}$
  for which $w$ belongs to $W'(d,2m;T(\pi))$.
Our immediate goal is to establish the following

\begin{theorem}\label{thm:ST}
There is a bijection between 
  $\barG_r(U,V;d)$ and the walks in $\widetilde{W}(d,2m;\vec{0})$
  staying in the region $x_1\geq x_2\geq \ldots \geq x_d$.
\end{theorem}
\medskip
We now discuss how to associate walks to Young tableaux.
First we need to introduce additional terminology.
We say that a walk $w=a_{1}\cdots a_{m}$ satisfies
\begin{quote}
\textbf{Condition (W):} If for each $i\in [n]$ and $1\leq s<r$ 
  it holds that $a_{r(i-1)+s} \geq a_{r(i-1)+s+1}$.
\end{quote}
Let $\varphi$ be the mapping from $T([m];d)$ to walks in $W(d,m;\ZZ^d)$
  such that $\varphi(T)=a_{1}\cdots a_{m}$ where $a_{i}$ equals the column
  in which entry $i$ appears in $T$.
It immediately follows that:
\begin{lemma}\label{lem:uno}
The mapping $\varphi$ is a bijection between tableaux in $T([m];d)$ satisfying 
  condition (T) and
  walks of length $m$ starting at the origin, moving only
  in positive directions, staying in the
  region $x_1\geq x_2\geq \dots \geq x_d$
  and satisfying condition (W).
\end{lemma}
\begin{proof}
If $\varphi(T)=\varphi(T')$ for $T,T'\in T([m];d)$, then $T$ and $T'$ 
  have the same elements in each of their columns. 
Since in a Young tableau the entries of each column are increasing
  from top to bottom, it follows that $T=T'$.
We have thus established that $\varphi$ is an injection.

Assume now that $w=a_{1}\cdots a_{m}$ is a walk 
  of length $m$ starting at the origin, moving only
  in positive directions, staying in the 
  region $x_1\geq x_2\geq \dots \geq x_d$ and satisfying condition (W).
Denote by $C(l)$ the set of indices $j$ for which $a_j=l$. 
Note that since $w$ is a walk in $\ZZ^d$, then $C(l)$ is empty for all $l>d$.
Let $T$ be the Young tableau whose $l$--th column entries correspond
  to $C(l)$ (obviously ordered increasingly from top to bottom).
Note that $T$ is indeed a Young tableau since $|C(1)|\geq |C(2)|\geq 
  \ldots \geq |C(d)|$ and given that the entries on each row of $T$ are
  strictly increasing (the latter follows from the fact that $w$ stays
  in the region $x_1\geq x_2\geq \dots \geq x_d$).
Observe that $T$ belongs to $T([m];d)$. 
We claim that $T$ satisfies condition (T).
Indeed, 
  by construction and since $w$ satisfies condition (W), for each $i\in[n]$
  it must hold that the indices of the columns of the 
  entries $r(i-1)+1,\ldots, ri$ of $T$ is a weakly decreasing sequence.
Hence, for every $1\leq s<r$, the entry $r(i-1)+s+1$ is weakly to the 
  left of $r(i-1)+s$.
Since $r(i-1)+s+1>r(i-1)+s$ and $T$ is a tableau, it must be the 
  case that the entry $r(i-1)+s+1$ is strictly below the entry
  $r(i-1)+s$.
\end{proof}
Note that if $T$ and $T'$ belong to $T([m];d)$ and have the same
  shape, then $\varphi(T)$ and $\varphi(T')$ are walks that terminate
  at the same lattice point.

\begin{corollary}\label{cor:dd}
There is a bijection between ordered pairs of tableaux
  of the same shape belonging to $T([m];d)$ satisfying condition (T),
  and walks in $\widetilde{W}(d,2m;\vec{0})$ staying in the region 
  $x_1\geq x_2\geq \ldots \geq x_d$.
\end{corollary}
\begin{proof}
By Lemma~\ref{lem:uno} there is a bijection between 
  ordered pairs of tableaux with the claimed properties and 
  ordered pairs of
  walks of length $m$ starting at the origin, moving only
  in positive directions that stay in the
  region $x_1\geq x_2\geq \ldots \geq x_d$
  and satisfy condition (W).
Say such pair of walks are $c_{1}\cdots c_{m}$ and 
  $c'_{1}\cdots c'_{m}$ respectively.
Then, $c_{1}\cdots c_{m}|c'_{m}\cdots c'_{1}$ 
  is the sought after walk with the desired properties.
\end{proof}
Figure~\ref{fig:young-walk} illustrates the bijection implicit in 
  the proof of Corollary~\ref{cor:dd}.
\begin{figure}
\begin{center}
\input young-walk.pstex_t
\end{center}
\caption{Walk in $\widetilde{W}(d,2m;\vec{0})$ associated to the 
  pair of Young tableaux of Figure~\ref{fig:perm-young} (and thus also
  to the graph of Figure~\ref{fig:example-graphs}).}
\label{fig:young-walk}
\end{figure}
Note that Theorem~\ref{thm:ST} is an immediate consequence of 
  Theorem~\ref{thm:perm-young} and Corollary~\ref{cor:dd}.

\begin{proof}{[of Theorem~\ref{thm.1}]}
The desired conclusion is an immediate consequence of 
  Theorem~\ref{thm:ST} and the existence of a 
  a parity-reversing involution $\rho$ on 
  the walks $w$ in 
  $\widetilde{W}(d,2m;\vec{0})$ not staying in the region 
  $x_1\geq x_2\geq \ldots \geq x_d$.
The involution is most easily described if we translate the walks to start at 
$(d-1, d-2,\dots, 0)$; the walks are then restricted not to lie completely
in the region $R$ defined by $x_1> x_2> \ldots > x_d$. 
Let $N$ be the subset of the translated walks of 
  $\widetilde{W}(d,2m;\vec{0})$ not lying completely in $R$.
Let $w= c_1 \ldots c_{2m}\in N$ and let $t$ be the smallest index such that 
  the walk given by the initial segment of $c_1\ldots c_t$ of $w$ 
  terminates in a vertex $(p_1,\ldots,p_d)\not\in R$. 
Hence, there is exactly one $j$ such that $p_j= p_{j+1}$.

Walk $\rho(w)$ is constructed as follows:
\begin{itemize}
\item Leave segment $c_1\ldots c_t$ unchanged.
\item For each $i\in [2n]$, define $S(i)=\set{s\in [2m]}{r(i-1)< s\leq ri}$,
  $S_{0}(i)=\set{s\in S(i)}{s>t, c_s=j}$ and
  $S_{1}(i)=\set{s\in S(i)}{s>t, c_s=j+1}$.
For $i\leq n$ (respectively $i>n$),
  assign the value $j+1$ to the $|S_0(i)|$ first (respectively last) 
  coordinates of $(c_{s} : s\in S_{0}(i)\cup S_{1}(i))$ and the value
  $j$ to the remaining $|S_{1}(i)|$ coordinates.
\end{itemize}
It is easy to see that if $w$ terminates in 
  $(q_1,\ldots, q_d)$, then $\rho(w)$ terminates in 
  $(q_1,\ldots, q_{j+1}, q_j, \ldots, q_{d})$. 
Hence, $\rho$ reverses the parity of $w$.
Moreover, $\rho\circ\rho$ is the identity.
It remains to show that $\rho(w)\in N$. 
Obviously $\rho(w)$ does not stay in $R$
  (as $w$ does not). 
Hence, it suffices to show the following: if 
  $\rho(w)= a_1\ldots a_m|b_1 \ldots b_m$, then for each 
  $i\in [n]$ and $1\leq s <r$ we have 
  $a_{r(i-1)+s}\geq a_{r(i-1)+s+1}$ and 
  $b_{r(i-1)+s}\leq b_{r(i-1)+s+1}$.
This is clearly true for every block $\set{r(i-1)+s}{1\leq s \leq r}$
  completely contained inside $w$'s unchanged segment (i.e., ${1,\ldots, t}$) 
  and inside $w$'s modified segment (i.e., ${t+1,\ldots,2m}$),
  given that it is true for $w$ and by the 
  definition of $\rho$.
There is still the case to handle where $t\in \set{r(i-1)+s}{1\leq s \leq r}$.
Here, it is true by the following observation: if $t\leq m$ then
  $c_t=j+1$, otherwise $c_t=j$.
\end{proof}

\section{Second proof}\label{sec:half-graphs}
Henceforth let $m=rn$. 
In this section we introduce two mappings $\Phi$ and $\phi$. 
The former is shown to be an injection that, when restricted to 
  $\mathcal{F}=\barG_r(U,V;d)$, takes values in $W'(d,2m;\vec{0})$.
Our first goal is to characterize those walks that belong to 
  $\Phi(\mathcal{F})$.
The second mapping $\phi$ plays a crucial role in fulfilling 
  this latter objective.
Then, relying on the aforementioned characterization we define a 
  parity reversing involution on $W'(d,2m;\Toep)\setminus\Phi(\mathcal{F})$.
This essentially establishes Theorem~\ref{thm.1}.

\medskip
Let $\Phi$ be the function that associates to an 
  $r$--configuration $F$ of $U$ and $V$ the value 
  $\Phi(F)= a_{u_1^1}\cdots a_{u_n^r}|b_{v_1^1}\cdots b_{v_n^r} 
  \in W^*(d,2m;\ZZ^d)$, where 
\begin{itemize}
\item
$a_{\baru}$ equals the largest size of a planar matching of 
  $F$ using nodes up to $\baru$,
\item
$b_{\barv}$ equals the largest size of a planar matching of 
  $F$ using nodes up to $\barv$.
\end{itemize}
Note that indeed $\Phi(F)\in W^*(d,2m;\ZZ^d)$ when $F$ is an
  $r$--configuration of $U$ and $V$. 
Figure~\ref{fig:walk} illustrates the definition of $\Phi(\cdot)$.
\begin{figure}
\begin{center}
\input walk.pstex_t
\end{center}
\caption{Walk $\Phi(\barG)=111122|112121$ for the multi--graph $G$ 
  of Figure~\ref{fig:example-graphs} 
  (direction $1$ is to the right and direction $2$ is up --- negative steps are represented by segmented lines.)}
\label{fig:walk}
\end{figure}

\begin{figure}
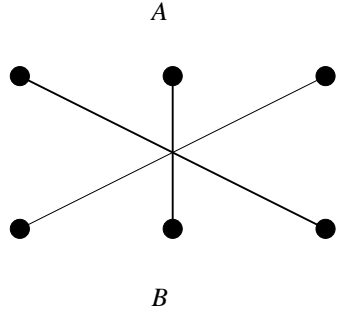

\begin{center}
\input crossing.pstex_t
\end{center}
\caption{Crossing obtained from left to right ordered sets $A$ and $B$.}
\label{fig:crossing}
\end{figure}

The following definition will be instrumental in the introduction
  of a mapping between walks and configurations.
\begin{definition}\label{def.cross}
Let $A$ and $B$ be two linearly ordered sets of equal size.
We say that a quasi configuration is obtained from 
  $A$ and $B$ \textit{in a crossing way} if the first element 
  of $A$ is paired with the last element of $B$, and so on, 
  until finally the last element of $A$ is paired to the first 
  element of $B$.
\end{definition}
Figure~\ref{fig:crossing} illustrates the concept just introduced.
We say that $H$ is a \textit{quasi $r$--configuration} of 
  $U$ and $V$ if it can be 
  obtained from a configuration $F$ of $U$ and $V$
  by ``breaking'' (deleting) some of its ``edges'' (pairings).
Note that the same quasi $r$--configuration may be obtained by ``breaking''
  different $r$--configurations.

\begin{figure}
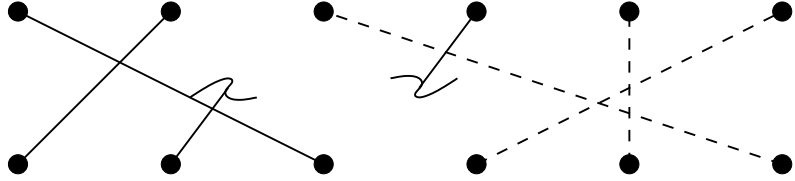

\begin{center}
\input phi.pstex_t
\end{center}
\caption{The quasi configuration $\phi(w)$ associated to $w=112122|122122$. 
Continuous lines corresponds to the crossing of $A_{1}(w)$ and 
  $B_{1}(w)$ and segmented lines to the crossing of $A_{2}(w)$
  and $B_{2}(w)$.}
\label{fig:phi}
\end{figure}

\medskip
For $w=a_{u_1^1}\cdots a_{u_n^r}|b_{v_1^1}\cdots b_{v_n^r}
  \in W^*(d,2m;\ZZ^d)$, let
  $A_k(w)=\set{\baru}{a_{\baru}= k}$ and 
  $B_k(w)=\set{\barv}{b_{\barv}= k}$.
We are now ready to introduce a mapping between walks and quasi configurations.
Let $\phi$ be a function that associates to a walk $w\in W^*(d,2m;\ZZ^d)$ 
  a quasi $r$--configuration $\phi(w)$ as follows:
  for each $k$, if  $|A_k(w)|\geq |B_k(w)|$ then connect the initial 
  segment of $A_k(w)$ of size $|B_k(w)|$ in a crossing way with 
  $B_k(w)$.
If  $|A_k(w)|\leq |B_k(w)|$, then connect the terminal segment of 
  $B_k(w)$ of size $|A_k(w)|$ in a crossing way with $A_k(w)$.
Figure~\ref{fig:phi} illustrates $\phi(\cdot)$'s definition.

\begin{fact}\label{fact:cross}
Let $w\in W^{*}(d,2m;\ZZ^d)$.
Then, the edges in $\phi(w)$ incident to two distinct elements of 
  $A_{k}(w)$ must cross.
A similar observation holds for $B_{k}(w)$.
\end{fact}

\begin{fact}\label{fact}
Let $w=\Phi(F)$ for an $r$--configuration $F$ of $U$ and $V$
  and let $(\baru,\barv)$ be a pairing of $F$.
Then, $\baru\in A_{k}(w)$ if and only if $\barv\in B_{k}(w)$.
\end{fact}

The following result gives an interpretation in terms of graphs
  of what it means for a walk starting at the origin to terminate
  also at the origin.
\begin{lemma}\label{lem:origin}
Let $w\in W^{*}(d,2m;\ZZ^d)$.
Then, $\phi(w)$ is an $r$--configuration of $U$ and $V$
  if and only if $w\in W^{*}(d,2m;\vec{0})$ for some $d$.
\end{lemma}
\begin{proof}
A closed walk $w$ passes through the origin if and only if
  $|A_{k}(w)|= |B_{k}(w)|$ for all $k$.
The latter is certainly equivalent to $\phi(w)$ being a configuration.
\end{proof}
We now prove a technical result.
\begin{lemma}\label{lem:tech}
Let $k\in [d]$ be arbitrary.
For every $r$--configuration $F$ of $U$ and $V$,
  the set of edges incident to $A_k(\Phi(F))$ equals 
  the set of edges incident to $B_k(\Phi(F))$.
\end{lemma}
\begin{proof}
Let $\baru\in A_k(\Phi(F))$. 
There is a unique $\barv$ such that $(\baru,\barv)$ is a pairing
  of $F$.
By Fact~\ref{fact}, it must hold that $\barv\in B_k(\Phi(F))$.  
\end{proof}
The following result establishes that $\Phi(\cdot)$ is an injection.
\begin{lemma}\label{lem:inj}
For every $r$--configuration $F$ of $U$ and $V$, 
  it holds that $\phi(\Phi(F))=F$.
\end{lemma}
\begin{proof}
By Fact~\ref{fact:cross} and Lemma~\ref{lem:tech}, $A_k(\Phi(F))$ 
  and  $B_k(\Phi(F))$ are equal size sets that 
  must be joined in the crossing way in $F$.
Since a pairing of $F$ is an element of $A_k(\Phi(F))\times B_k(\Phi(F))$
  for some $k$, it follows that $\phi(\Phi(F))=F$.
\end{proof}

\begin{lemma}\label{lem:2}
Let $\mathcal{F}$ be a family of $r$--configurations of $U$ and $V$.
A walk $w$ belongs to 
  $\Phi(\mathcal{F})$ if and only if $\Phi(\phi(w))= w$ and 
  $\phi(w)\in \mathcal{F}$.
\end{lemma}
\begin{proof}
If $\phi(w)\in \mathcal{F}$, then $w=\Phi(\phi(w))$
  belongs to $\Phi(\mathcal{F})$.
If  $w= \Phi(F)$ for some $r$--configuration $F$ of $U$ and $V$, then
  Lemma~\ref{lem:inj} implies that $\phi(w)=F$. 
If in addition $F\in \mathcal{F}$, then 
  one gets that $\phi(w)\in \mathcal{F}$.
\end{proof}
Two walks in $W^{*}(d,2m;\ZZ^d)$ are certainly equal if their sequence of
  positive and negative steps agree.
The next lemma gives a simpler 
  necessary and sufficient condition for the equality of 
  two walks $w$ and $\Phi(\phi(w))$ when $w$ is a closed walk that goes
  through the origin.
Indeed, it says that one only needs to focus on establishing the 
  equality of the sequence of their positive steps.
The result will be useful later in order to establish the equality of 
  two walks $w$ and $\Phi(\phi(w))$.
\begin{lemma}\label{lem:eq-positive}
Let $w\in W^{*}(d,2m;\vec{0})$.
Then, $\Phi(\phi(w))$ and $w$ agree in their positive steps
  if and only if $\Phi(\phi(w))=w$.
\end{lemma}
\begin{proof}
If  $\Phi(\phi(w))=w$, then $\Phi(\phi(w))$ and $w$ clearly 
  agree in their positive steps.
To prove the converse, let $w'=\Phi(\phi(w))$.
Assume $w'$ and $w$ agree in their positive steps.
First, recall that by Lemma~\ref{lem:origin}, $\phi(w)$ is an
  $r$--configuration of $U$ and $V$.
Hence, Lemma~\ref{lem:inj} implies that $\phi(w')=\phi(\Phi(\phi(w)))=\phi(w)$.
Thus, $A_k(w) = A_k(w')$ for every $k$. 
Since $\phi(w')$ and $\phi(w)$ are the same configurations,  
  they have the same set of edges.
Consider $\barv\in B_{k}(w)$.
There is a unique edge $(\baru,\barv)$ of $\phi(w)$ incident on $\barv$.
By Fact~\ref{fact}, we have that $\baru\in A_{k}(w)=A_{k}(w')$.
But edge $(\baru,\barv)$ is an edge of $\phi(w')$.
Hence, again by Fact~\ref{fact}, we get that $\barv\in B_{k}(w')$.
We have shown that $B_{k}(w)\subseteq B_{k}(w')$.
The reverse inclusion can be similarly proved.
Since $k$ was arbitrary, we conclude that the negative
  steps of $w$ and $w'$ are the same,
  and the two walks must thus be equal.
\end{proof}

For the walk $w=a_{u_1^1}\cdots a_{u_n^r}|b_{v_1^1} \cdots b_{v_n^r}$, 
  denote by $k(\baru)$ and $l(\baru)$ the number of occurrences of 
  $a_{\baru}$ and $a_{\baru}-1$ in 
  $\set{\baru'\in\barU}{\baru'\preceq \baru}$ respectively.
\begin{example}
For the walk $111122|112121$ of Figure~\ref{fig:walk} we have:
\begin{center}
\begin{tabular}{|c||rrrrrr|}\hline
$\baru$  & $1$ & $2$ & $3$ & $4$ & $5$ & $6$ \\ \hline
$k(\baru)$ & $1$ & $2$ & $3$ & $4$ & $1$ & $2$ \\
$l(\baru)$ & $0$ & $0$ & $0$ & $0$ & $4$ & $4$ \\ \hline
\end{tabular}
\end{center}
\end{example}
We say that $w$ satisfies 
\begin{quote}
\textbf{Condition (C):} If
  for each $\baru$ such that 
  $a_{\baru}> 1$, $l(\baru)>0$ and the $l(\baru)$--th-to-last appearance of 
  $a_{\baru}-1$ in the negative steps of $w$, if it exists, comes 
  before the $k(\baru)$--th-to-last appearance of $a_{\baru}$ in the negative 
  steps of~$w$. 
\end{quote}

\begin{lemma}\label{lem:condition}
Let $w= a_{u_1^1}\cdots a_{u_n^r}|b_{u_1^r}\cdots b_{v_n^r}$
  be a walk in $W^{*}(d,2m;\Toep)$.
Then, $\Phi(\phi(w))= w$ and $\phi(w)$ is an $r$--configuration of $U$ and $V$
  if and only if $w$ satisfies condition (C).
\end{lemma}
\begin{proof}
Since
  $\phi(w)$ is an $r$--configuration of $U$ and $V$,
  Lemma~\ref{lem:origin} implies that $w$ must terminate at the origin.
On the other hand, if $w$ satisfies condition (C),
  then $w$ also needs to terminate at the origin.
Hence, $\phi(w)$ would be an $r$--configuration of $U$ and $V$.
Indeed, let $j$ be 
  the smallest coordinate in which the terminal point of $w$ is positive.
Note that $j> 1$ by the definition of Toeplitz points. 
Let $\baru$ be maximum so that $a_{\baru}= j$. 
By the choice of $j$, there will be fewer than $k(\baru)$ appearances 
  of $j$ and at least $l(\baru)$ appearances  of $j-1$ among the negative 
  steps.
This contradicts the the fact that $w$ satisfies condition (C).
We thus can assume without loss of generality that $w$ starts and 
  ends at the origin.

Let $\Phi(\phi(w)) = a'_{u_1^1}\cdots a'_{u_n^r}|b'_{v_1^r}\cdots b'_{v_n^r}$.
By Lemma~\ref{lem:eq-positive}, 
  $\Phi(\phi(w))$ and $w$ are distinct if and only if
  they differ in some positive step.
Assume that $a_{\baru}= a'_{\baru}$ 
  for each $\baru\prec\tilde{u}$ and 
  $a_{\tilde{u}}\neq a'_{\tilde{u}}$. 
We claim that $a'_{\tilde{u}}\leq a_{\tilde{u}}$.
Indeed, suppose this is not the case.
Since $a'_{\tilde{u}}$ equals the size of the largest  
  planar matching of $\phi(w)$ up to $\tilde{u}$,
  there is a $\baru\prec\tilde{u}$ such that
  $a'_{\baru}=a_{\tilde{u}}$ and the
  edges incident to $\baru$ and 
  $\tilde{u}$ of $\phi(w)$ are non-crossing.  
We also have $a'_{\baru}= a_{\baru}$ by the choice of $\tilde{u}$.
Hence, $\baru$ and $\tilde{u}$ belong to $A_k(w)$ for some $k$.
By Fact~\ref{fact:cross} the edges incident to  
  $\baru$ and $\tilde{u}$ must be non-crossing.
A contradiction.
This establishes our claim.

It follows that
  $a'_{\tilde{u}}= a_{\tilde{u}}$ if and 
  only if $a'_{\tilde{u}}\geq a_{\tilde{u}}$.
We now establish a condition equivalent to 
  $a'_{\tilde{u}}\geq a_{\tilde{u}}$
  by considering the following two cases:
\begin{itemize}
\item \textbf{Case $a_{\tilde{u}}= 1$:}
Then, certainly $a'_{\tilde{u}}\geq a_{\tilde{u}}$.

\item \textbf{Case $a_{\tilde{u}}>1$:} Then, there 
  is a $\baru\prec\tilde{u}$ such that 
  $a'_{\baru}= a_{\tilde{u}}-1$ and the edges incident to
  $\baru$ and $\tilde{u}$ are non--crossing in $\phi(w)$. 
So we can extend with the edge incident to $\tilde{u}$ 
  the size $a'_{\baru}$ planar matching of $\phi(w)$
  up to $\baru$.
Thus, it must be the case that 
  $a'_{\tilde{u}}\geq a_{\tilde{u}}$.
\end{itemize}
Summarizing $a_{\tilde{u}}=a'_{\tilde{u}}$ if 
  and only if
\begin{itemize}
\item
$a_{\tilde{u}}=1$, or
\item 
if $a_{\tilde{u}}>1$ and there is a $\baru\prec\tilde{u}$ such that 
  $a'_{\baru}= a_{\tilde{u}}-1$ and the edges 
  $\tilde{u}$ and $\baru$ are non--crossing in $\phi(w)$.
\end{itemize}
The lemma follows by observing that when 
  $a_{\tilde{u}}>1$, 
  the fact that $w$ satisfies condition (C) amounts 
  to saying that there is a $\baru\prec \tilde{u}$ such that
  $a_{\baru}=a_{\tilde{u}}-1$ and the edges incident
  to $\baru$ and $\tilde{u}$ are non--crossing in $\phi(w)$.
So, all positive steps of $\Phi(\phi(w))$ and $w$ agree 
  if and only if 
  for each $\baru$ such that 
  $a_{\baru}> 1$, $l(\baru)>0$ and the $l(\baru)$--th-to-last appearance of 
  $a_{\baru}-1$ in the negative steps of $w$, if it exists, comes 
  before the $k(\baru)$--th-to-last appearance of $a_{\baru}$ in the negative 
  steps of $w$. 
\end{proof}

So far in this section we have not directly being concerned with 
  walks $W'(d,2m;\Toep)$ nor the collection of configurations 
  $\barG_r(U,V;d)$.
The next result is the link through which we use all previous 
  results in order to prove Theorem~\ref{thm.1}.
\begin{lemma}\label{lem:link}
Let $w\in W'(d,2m;\vec{0})$. 
If $\Phi(\phi(w))=w$, then $\phi(w)\in\barG_r(U,V;d)$.
\end{lemma}
\begin{proof}
Suppose $\Phi(\phi(w))=w$ and 
  $w=a_{u_1^1}\cdots a_{u_n^r}|b_{v_1^1}\cdots b_{v_n^r}$
  is such that $\phi(w)$ does not belong to $\barG_r(U,V;d)$.
Note that since $w$ is a closed walk that goes through the 
  origin, by Lemma~\ref{lem:origin}, we have that $\phi(w)$ 
  is an $r$--configuration of $U$ and $V$.
Thus, it must be the case that 
  either there is a $u\in U$ such that for some $s< t$ the edges
  incident to $u^s$ and $u^t$ are non--crossing, or
  there  is a $v\in V$ such that for some $s< t$ the edges
  incident to $v^s$ and $v^t$ are non--crossing.
Without loss of generality assume the former case holds.
It follows that, the largest planar matching up to $u^s$ is strictly smaller
  than the largest planar matching up to $u^t$, i.e.,
  $a_{u^s} < a_{u^t}$.
This contradicts the fact that $w$ belongs to $W'(d,2m;\vec{0})$.
\end{proof}

\begin{theorem}\label{thm:bijec}
The mapping $\Phi$ is a bijection between $\barG_r(U,V;d)$
  and the collection of walks in $W'(d,2m;\Toep)$ satisfying
  condition (C).
\end{theorem}
\begin{proof}
By Lemma~\ref{lem:inj} we know that $\Phi$ is an injection.
We claim it is also onto. 
Indeed, if $w$ is a walk in $W'(d,2m;\Toep)$ satisfying condition (C),
  then Lemmas~\ref{lem:origin}, 
  \ref{lem:2}, \ref{lem:condition}, and \ref{lem:link}
  imply that $\phi(w)$ belongs to $\barG_r(U,V;d)$ and
  $\Phi(\phi(w))=w$.  
\end{proof}

\begin{proof}{[of Theorem~\ref{thm.1}]}
The desired conclusion is an immediate consequence of 
  Theorem~\ref{thm:bijec} and the existence of 
  a parity-reversing involution $\rho$ on 
  walks $w$ in $W'(d,2m;\Toep)$ that don't satisfy condition (C).
To define $\rho$, assume
  $w= a_{u_1^1}\cdots a_{u_n^r}|b_{v_1^1}\cdots b_{v_n^r}$ 
  and let $\baru$ be the smallest index for which $w$ does not satisfy
  condition (C).
Let $\barv$ be such that $b_{\barv}$ is the 
  $l(\baru)$--th-to-last occurrence of $a_{\baru}-1$ among the 
  negative steps; if $l(\baru) = 0$ then let $\barv= rn+1$.

Walk $\rho(w)$ is constructed as follows:
\begin{itemize}
\item Leave segments $a_{u_1^1}\cdots a_{\baru}$ and 
  $b_{\barv} \cdots b_{v_n^r}$ unchanged.
\item For every $i\in[n]$, let  
  $S_0(i) = \set{s}{a_{u_{i}^{s}}=a_{\baru},\baru\prec u_{i}^{s}}$ and
  $S_{1}(i) = \set{s}{a_{u_{i}^{s}}=a_{\baru}{-}1,\baru\prec u_{i}^{s}}$.
Assign the value $a_{\baru}$ to the $|S_{1}(i)|$ first coordinates in 
  $(a_{u_{i}^{s}} : s\in S_{0}(i)\cup S_{1}(i))$ and 
  the value $a_{\baru}-1$ to the remaining $|S_{0}(i)|$ coordinates.
\end{itemize}
The application of $\rho$ does not change the smallest index 
  not satisfying the sufficient condition of Lemma~\ref{lem:condition}.
It follows that $\rho(w)$ also violates condition (C).
We claim that $\rho(w)=a'_{u_1^1}\cdots a'_{u_n^r}|b'_{v_1^1}\cdots b'_{v_n^r}$
  belongs to $W'(d,2m;\Toep)$.
We need to show that for each 
  $i\in[n]$ and $1\leq s <r$, we have 
  $a'_{u_i^{s}} \geq a'_{u_i^{s+1}}$ 
  and $b'_{v_i^{s}}\geq b'_{v_i^{s+1}}$.
This is clearly true for every block 
  $\set{u_i^s}{s\in[r]}$ completely contained in the unchanged segments, 
  and also inside the modified segment.
The remaining two cases to consider are $\baru=u_i^s$ and/or $\barv=v_j^{s'}$
  for some $s<r$ and/or $s'>1$.
Both cases are easy to handle.
We leave the details to the reader.
 
Assume $w$ terminates at $T(\pi)$ for some permutation $\pi$ of $[d]$.
Let $\tau$ be a transposition of $a_{\baru}$ and $a_{\baru}-1$.  
Finally, we claim that $\rho(w)$ terminates in $T(\pi\circ\tau)$.
Indeed, by our choice of $\baru$, the number of appearances of 
  $a_{\baru}$ in $b_{\barv}\cdots b_{v_n^r}$ is less
  than $k(\baru)$. 
It must equal to $k(\baru)-1$, otherwise we could have chosen the index
  of the $(k(\baru)-1)$--th appearance of $a_{\baru}$ for $\baru$. 
Hence, in the unchanged segments of the walk $w$,
  there is one net positive step in direction $a_{\baru}$ and zero 
  net steps in direction $a_{\baru}-1$. 
It follows that in the segment of $w$ that changes, there 
  are $a_{\baru}-\pi(a_{\baru})-1$  and 
  $a_{\baru}-1-\pi(a_{\baru}-1)$ net positive steps in directions
  $a_{\baru}$ and $a_{\baru}-1$  respectively.
Let $\sigma_s$ denote the $s$--th coordinate of 
  the terminal point of a walk $\sigma$. 
We get that $\rho(w)_{a_{\baru}}
  = a_{\baru}-1-\pi(a_{\baru}-1) + 1= a_{\baru}- \pi(a_{\baru}-1)$
  and similarly $\rho(w)_{a_{\baru}-1}
  = a_{\baru} - \pi(a_{\baru}) - 1= a_{\baru}- 1- \pi(a_{\baru})$.
\end{proof}

\bibliographystyle{alpha}
\bibliography{biblio}

\end{document}